\overfullrule=0pt
\centerline {\bf A class of functionals possessing multiple global minima}\par
\bigskip
\bigskip
\centerline {BIAGIO RICCERI}\par
\bigskip
\bigskip
\centerline {\it To Professor Gheorghe Morosanu, with friendship, on his 70th birthday}\par
\bigskip
\bigskip
{\bf Abstract.}We get a new multiplicity result for gradient systems. Here is a very particular corollary: Let $\Omega\subset
{\bf R}^n$ ($n\geq 2$) be a smooth bounded domain and let $\Phi:{\bf R}^2\to {\bf R}$ be a $C^1$ function, with $\Phi(0,0)=0$,
such that 
$$\sup_{(u,v)\in {\bf R}^2}{{|\Phi_u(u,v)|+|\Phi_v(u,v)|}\over {1+|u|^p+|v|^p}}<+\infty$$
where $p>0$, with $p={{2}\over {n-2}}$ when $n>2$.\par
Then, for every convex set $S\subseteq L^{\infty}(\Omega)\times L^{\infty}(\Omega)$ dense in $L^2(\Omega)\times
L^2(\Omega)$, there exists $(\alpha,\beta)\in S$ such that the problem
 $$\cases {-\Delta u=(\alpha(x)\cos(\Phi(u,v))-\beta(x)\sin(\Phi(u,v)))\Phi_u(u,v) & in $\Omega$
\cr & \cr
-\Delta v= (\alpha(x)\cos(\Phi(u,v))-\beta(x)\sin(\Phi(u,v)))\Phi_v(u,v) & in $\Omega$
\cr & \cr
u=v=0 & on $\partial\Omega$\cr}$$
has at least three weak solutions, two of which are global minima in $H^1_0(\Omega)\times H^1_0(\Omega)$ of the functional
$$(u,v)\to {{1}\over {2}}\left ( \int_{\Omega}|\nabla u(x)|^2dx+\int_{\Omega}|\nabla v(x)|^2dx\right )$$
$$-\int_{\Omega}(\alpha(x)\sin(\Phi(u(x),v(x)))+\beta(x)\cos(\Phi(u(x),v(x))))dx\ .$$
\bigskip
\bigskip
\bigskip
\bigskip
{\bf 1. Introduction}\par
\bigskip
Let $S$ be a topological space. A function $g:S\to {\bf R}$ is said to be inf-compact if, for each $r\in {\bf R}$, the set
$g^{-1}(]-\infty,r])$ is compact.\par
\smallskip
If $Y$ is a real interval and $f:S\times Y\to {\bf R}$ is a function inf-compact and lower semicontinuous in $S$, and concave in $Y$,
the occurrence of the strict minimax inequality
$$\sup_Y\inf_Sf<\inf_S\sup_Yf$$
implies that the interior of the set $A$ of all $y\in Y$ for which $f(\cdot,y)$ has at least two local minima is non-empty.
This fact was essentially shown in [4], giving then raise to an enormous number of subsequent applications to the multiplicity
of solutions for nonlinear equations of variational nature (see [7] for an account up to 2010).\par
\smallskip
In [6] (see also [5]), we realized that, under the same assumptions as above, the occurrence of the strict minimax inequality also implies the existence
of $\tilde y\in Y$ such that the function $f(\cdot,\tilde y)$ has at least two global minima. It may happen that $\tilde y$ is unique and does not belong
to the closure of $A$ (see Example 7 of [1]).\par
\smallskip
In [8] and [12], we extended the result of [6] to the case where $Y$ is an arbitrary convex set in a vector space. We also stress that such an extension is not possible for the result of [4]. We then started to build a network of applications of the results of [8] and [12] which touches several different topics: uniquely remotal sets in normed spaces ([8]); non-expansive operators ([9]); singular points ([10]); Kirchhoff-type problems ([11]); 
Lagrangian systems of relativistic oscillators ([13]); integral functional of the Calculus of Variations ([14]);  non-cooperative gradient systems ([15]); variational inequalities ([16]).
\par
\smallskip
The aim of this short paper is to establish a further application within that network.\par
\bigskip
{\bf 2. Results}\par
\bigskip
The main abstract result is as follows:\par
\medskip
THEOREM 2.1. - {\it Let $X$ be a topological space, $(Y,\langle\cdot,\cdot,\rangle)$ a real Hilbert space, $T\subseteq Y$ a convex
set dense in $Y$ and $I:X\to {\bf R}$,
$\varphi:X\to Y$ two functions such that, for each $y\in T$, the function $x\to I(x)+\langle \varphi(x),y\rangle$ is lower semicontinuous
and inf-compact. Moreover, assume that there exists a point $x_0\in X$, with $\varphi(x_0)\neq 0$, such that\par
\noindent
$(a)$\hskip 5pt $x_0$ is a global minimum of both functions $I$ and $\|\varphi(\cdot)\|$\ ;\par
\noindent
$(b)$\hskip 5pt $\inf_{x\in X}\langle\varphi(x),\varphi(x_0)\rangle<\|\varphi(x_0)\|^2$\ .\par
Then, for each convex set $S\subseteq T$ dense in $Y$, there exists $y^*\in S$ such that the function $x\to I(x)+\langle\varphi(x),y^*\rangle$ has at least two global minima in $X$.}\par
\smallskip
PROOF. In view of $(b)$, we can find $\tilde x\in X$ and $r>0$ such that
$$I(\tilde x)+{{r}\over {\|\varphi(x_0)\|}}\langle \varphi(\tilde x),\varphi(x_0)\rangle<I(x_0)+r\|\varphi(x_0)\|\ .\eqno{(1)}$$
Thanks to $(a)$, we have
$$I(x_0)+r\|\varphi(x_0)\|=\inf_{x\in X}(I(x)+r\|\varphi(x)\|)\ .\eqno{(2)}$$
The function $y\to \inf_{x\in X}(I(x)+\langle\varphi(x),y\rangle)$ is weakly upper semicontinuous, and so there exists
$\tilde y\in B_r$ such that
$$\inf_{x\in X}(I(x)+\langle\varphi(x),\tilde y\rangle)=\sup_{y\in B_r}\inf_{x\in X}(I(x)+\langle\varphi(x),y\rangle)\ ,\eqno({3)}$$
$B_r$ being the closed ball in $X$, centered at $0$, of radius $r$.
We distinguish two cases. First, assume that $\tilde y\neq {{r\varphi(x_0)}\over {\|\varphi(x_0\|}}$. As a consequence,
taking into account that $r\|\varphi(x_0)\|$ is the maximum of the restriction to $B_r$ of the continuous linear
functional $\langle\varphi(x_0),\cdot\rangle$ (attained at the point ${{r\varphi(x_0)}\over {\|\varphi(x_0)\|}}$ only),
we have
$$\inf_{x\in X}(I(x)+\langle\varphi(x),\tilde y\rangle)\leq I(x_0)+\langle\varphi(x_0),\tilde y\rangle<I(x_0)+r\|\varphi(x_0)\|\ .\eqno{(4)}$$
Now, assume that $\tilde y={{r\varphi(x_0)}\over {\|\varphi(x_0\|}}$. In this case, due to $(1)$, we have
$$\inf_{x\in X}(I(x)+\langle\varphi(x),\tilde y\rangle)\leq I(\tilde x)+\langle \varphi(\tilde x),\tilde y\rangle =
I(\tilde x)+{{r}\over {\|\varphi(x_0)\|}}\langle \varphi(\tilde x),\varphi(x_0)\rangle<I(x_0)+r\|\varphi(x_0)\|\ .\eqno{(5)}$$
Therefore, from $(2)$, $(3)$, $(4)$ and $(5)$, it follows that
$$\sup_{y\in B_r}\inf_{x\in X}(I(x)+\langle\varphi(x),y\rangle)<\inf_{x\in X}\sup_{y\in B_r}(I(x)+\langle\varphi(x),y\rangle)\ .\eqno{(6)}$$
Now, let $S\subseteq T$ be a convex set dense in $Y$. By continuity, we clearly have
$$\sup_{y\in B_r\cap S}\langle\varphi(x),y\rangle=\sup_{y\in B_r}\langle\varphi(x),y\rangle$$
for all $x\in X$. Therefore, in view of $(6)$, we have
$$\sup_{y\in B_r\cap S}\inf_{x\in X}(I(x)+\langle\varphi(x),y\rangle)\leq
\sup_{y\in B_r}\inf_{x\in X}(I(x)+\langle\varphi(x),y\rangle)<\inf_{x\in X}\sup_{y\in B_r}(I(x)+\langle\varphi(x),y\rangle)=
\inf_{x\in X}\sup_{y\in B_r\cap S}(I(x)+\langle\varphi(x),y\rangle)\ .$$
At this point, the conclusion follows directly applying Theorem 1.1 of [12] to the restriction of the function
$(x,y)\to I(x)+\langle\varphi(x),y\rangle$ to $X\times (B_r\cap S)$.\hfill $\bigtriangleup$\par
\medskip
We now present an application of Theorem 2.1 to elliptic systems.\par
\smallskip
In the sequel, $\Omega\subseteq {\bf R}^n$ ($n\geq 2$) is a bounded domain with smooth boundary.\par
\smallskip
We denote by ${\cal A}$ the class of all functions $H:\Omega\times {\bf R}^2\to {\bf R}$ which are measurable in $\Omega$,
$C^1$ in ${\bf R}^2$ and satisfy
$$\sup_{(x,u,v)\in \Omega\times {\bf R}^2}{{|H_u(x,u,v)|+|H_v(x,u,v)|}\over {1+|u|^p+|v|^p}}<+\infty$$
where $p>0$, with $p<{{n+2}\over {n-2}}$ when $n>2$. \par
\smallskip
Given $H\in {\cal A}$, we are interested in the problem\par
$$\cases {-\Delta u=H_u(x,u,v) & in $\Omega$
\cr & \cr
-\Delta v=H_v(x,u,v) & in $\Omega$
\cr & \cr
u=v=0 & on $\partial\Omega$\ ,\cr} \eqno{(P_H)}$$
$H_u$ (resp. $H_v$) denoting the derivative of $H$ with respect to $u$ (resp. $v$).\par
\smallskip
As usual, a weak solution of $(P_H)$ is any $(u,v)\in H^1_0(\Omega)\times H^1_0(\Omega)$ such that
$$\int_{\Omega}\nabla u(x)\nabla\varphi(x)dx=\int_{\Omega}H_u(x,u(x),v(x))\varphi(x)dx\ ,$$
$$\int_{\Omega}\nabla v(x)\nabla\psi(x)dx=\int_{\Omega}H_v(x,u(x),v(x))\psi(x)dx$$
for all $\varphi, \psi\in H^1_0(\Omega)$.\par
\smallskip
Define the functional $I_H:H^1_0(\Omega)\times H^1_0(\Omega)\to {\bf R}$ by
$$I_H(u,v)={{1}\over {2}}\left ( \int_{\Omega}|\nabla u(x)|^2dx+\int_{\Omega}|\nabla v(x)|^2dx\right )-
\int_{\Omega}H(x,u(x),v(x))dx$$
for all $(u,v)\in H^1_0(\Omega)\times H^1_0(\Omega)$.\par
\smallskip
Since $H\in {\cal A}$, the functional $I_H$ is $C^1$ in $H^1_0(\Omega)\times H^1_0(\Omega)$ and its critical
points are precisely the weak solutions of $(P_H)$.  Moreover, due to the Sobolev embedding theorem, the functional
$(u,v)\to \int_{\Omega}H(x,u(x),v(x))$ has a compact derivative and, as a consequence, it is sequentially weakly continuous
in $H^1_0(\Omega)\times H^1_0(\Omega)$.\par
\smallskip
Also, we denote by $\lambda_1$ the first eigenvalue of the Dirichlet problem
$$\cases{-\Delta u=\lambda u & in $\Omega$ \cr & \cr u=0 & on $\partial\Omega$\ .\cr}$$
\smallskip
Our result is as follows:\par
\medskip
THEOREM 2.2. - {\it Let $F, G\in {\cal A}$, with $p={{2}\over {n-2}}$ when $n>2$, and let $K\in {\cal A}$,
with $K(x,0,0)=0$ for all $x\in \Omega$, satisfy the following conditions:\par
\noindent
$(a_1)$\hskip 5pt one has
$$\lim_{s^2+t^2\to +\infty}{{\sup_{x\in \Omega}(|F(x,s,t)|+|G(x,s,t)|)}\over {s^2+t^2}}=0\ ;$$
\noindent
$(a_2)$\hskip 5pt there is $\eta\in \left ]0,{{\lambda_1}\over {2}}\right [$ such that
$$K(x,s,t)\leq \eta(s^2+t^2)$$
for all $x\in \Omega$, $s, t\in {\bf R}$\ ;\par
\noindent
$(a_3)$\hskip 5pt one has
$$\hbox {\rm meas}(\{ x\in\Omega: 0<|F(x,0,0)|^2+|G(x,0,0)|^2\})>0 \eqno{(7)}$$
and
$$|F(x,0,0)|^2+|G(x,0,0)|^2\leq |F(x,s,t)|^2+|G(x,s,t)|^2 \eqno{(8)}$$
for all $x\in \Omega$, $s, t\in {\bf R}$\ ;\par
\noindent
$(a_4)$\hskip 5pt one has $$\hbox {\rm meas}\left (\left \{x\in \Omega : \inf_{(s,t)\in {\bf R}^2}(F(x,0,0)F(x,s,t)+G(x,0,0)G(x,s,t))<
|F(x,0,0)|^2+|G(x,0,0)|^2\right \}\right )>0\ .$$
Then, for every convex set $S\subseteq L^{\infty}(\Omega)\times L^{\infty}(\Omega)$ dense in $L^2(\Omega)\times
L^2(\Omega)$, there exists $(\alpha,\beta)\in S$ such that the problem
$$\cases {-\Delta u=\alpha(x)F_u(x,u,v)+\beta(x)G_u(x,u,v)+K_u(x,u,v) & in $\Omega$
\cr & \cr
-\Delta v= \alpha(x)F_v(x,u,v)+\beta(x)G_v(x,u,v)+K_v(x,u,v) & in $\Omega$
\cr & \cr
u=v=0 & on $\partial\Omega$\cr}$$
has at least three weak solutions, two of which are global minima in $H^1_0(\Omega)\times H^1_0(\Omega)$ of the functional
$$(u,v)\to {{1}\over {2}}\left ( \int_{\Omega}|\nabla u(x)|^2dx+\int_{\Omega}|\nabla v(x)|^2dx\right )-
\int_{\Omega}(\alpha(x)F(x,u(x),v(x))+\beta(x)G(x,u(x),v(x))+K(x,u(x),v(x)))dx\ .$$}\par
\smallskip
PROOF. We are going to apply Theorem 2.1, with the following choices: $X$ is the space $H^1_0(\Omega)\times H^1_0(\Omega)$ endowed
with the weak topology induced by the scalar product
$$\langle (u,v), (w,\omega)\rangle_X=\int_{\Omega}(\nabla u(x)\nabla w(x)+\nabla v(x)\nabla\omega(x))dx\ ;$$
$Y$ is the space $L^2(\Omega)\times L^2(\Omega)$ with the scalar product
$$\langle (f,g), (h,k)\rangle_Y=\int_{\Omega}(f(x)h(x)+g(x)k(x))dx\ ;$$
$T$ is $L^{\infty}(\Omega)\times L^{\infty}(\Omega)$; $I$ is the function defined by
$$I(u,v)={{1}\over {2}}\left ( \int_{\Omega}|\nabla u(x)|^2dx+\int_{\Omega}|\nabla v(x)|^2dx\right )-\int_{\Omega}K(x,u(x),v(x))dx$$
for all $(u,v)\in X$; $\varphi$ is the function defined by
$$\varphi(u,v)=(F(\cdot,u(\cdot),v(\cdot)),G(\cdot,u(\cdot),v(\cdot)))$$
for all $(u,v)\in X$; $x_0$ is the zero of $X$. Let us show that the assumptions of
Theorem 2.1 are satisfied. First, from $(7)$ and $(8)$ it clearly follows, respectively, that
$$\|\varphi(0,0)\|_Y^2=\int_{\Omega}(|F(x,0,0)|^2+|G(x,0,0|^2)dx>0$$
and that
$$\|\varphi(0,0)\|_Y^2\leq \|\varphi(u,v)\|_Y^2$$
for all $(u,v)\in X$. Moreover, from $(a_2)$, thanks to the Poincar\'e inequality, we get
$$\int_{\Omega}K(x,u(x),v(x))dx\leq \eta\int_{\Omega}(|u(x)|^2+|v(x)|^2)dx\leq
{{\eta}\over {\lambda_1}}\int_{\Omega}(|\nabla u(x)|^2+|\nabla v(x)|^2)dx \eqno{(9)}$$
for all $(u,v)\in X$. In particular, since $K(x,0,0)=0$ in $\Omega$ and ${{\eta}\over {\lambda_1}}<{{1}\over {2}}$,
from $(9)$ we infer that $(0,0)$ is a global minimum of $I$
in $X$. So, condition $(a)$ is satisfied. Now, let us verify condition $(b)$. To this end, set
$$P(x,s,t)=F(x,0,0)F(x,s,t)+G(x,0,0)G(x,s,t)-|F(x,0,0)|^2-|G(x,0,0)|^2$$
for all $(x,s,t)\in \Omega\times {\bf R}^2$ and
$$D=\left\{x\in \Omega : \inf_{(s,t)\in {\bf R}^2}P(x,s,t)<0\right\}\ .$$
By $(a_4)$, $D$ has a positive measure. In view of the Scorza-Dragoni theorem, there exists a compact set $C\subset D$, with
positive measure, such that the restriction of $P$ to $C\times {\bf R}^2$ is continuous. Fix a point $\tilde x\in C$ such that
the intersection of $C$ and any ball centered at $\tilde x$ has a positive measure. Choose $\tilde s,\tilde t\in {\bf R}\setminus \{0\}$
so that $P(\tilde x,\tilde s,\tilde t)<0$. By continuity, there is $r>0$ such that
$$P(x,\tilde s,\tilde t)<0$$
for all $x\in C\cap B(\tilde x,r)$. Set
$$\gamma=\sup_{(x,s,t)\in \Omega\times [-|\tilde s|,|\tilde s|]\times [-|\tilde t|,|\tilde t|]}|P(x,t,s)|\ .$$
Since $F, G\in {\cal A}$, $\gamma$ is finite. Now, choose an open set $A$ such that
$$C\cap B(\tilde x,r)\subset A\subset \Omega$$
and
$$\hbox {\rm meas}(A\setminus (C\cap B(\tilde x,r)))<-{{\int_{C\cap B(\tilde x,r)}P(x,\tilde s, \tilde t)dx}\over {\gamma}}\ .\eqno{(10)}$$
Finally, choose two functions $\tilde u, \tilde v\in H^1_0(\Omega)$ such that
$$\tilde u(x)=\tilde s\ ,\hskip 5pt \tilde v(x)=\tilde t$$
for all $x\in C\cap B(\tilde x,r)$\ ,\par
$$\tilde u(x)=\tilde v(x)=0$$
for all $x\in \Omega\setminus A$ and
$$|\tilde u(x)|\leq |\tilde s|\ ,\hskip 5pt |\tilde v(x)|\leq |\tilde t|$$
for all $x\in \Omega$. Then, taking $(10)$ into account, we have
$$\langle\varphi(\tilde u,\tilde v),\varphi(0,0)\rangle_Y-\|\varphi(0,0)\|_Y^2=\int_{\Omega}P(x,\tilde u(x),\tilde v(x))dx=
\int_{C\cap B(\tilde x,r)}P(x,\tilde s,\tilde t)dx+\int_{A\setminus (C\cap B(\tilde x,r))}P(x,\tilde u(x),\tilde v(x))dx$$
$$<\int_{C\cap B(\tilde x,r)}P(x,\tilde s,\tilde t)dx+\gamma\hbox {\rm meas}(A\setminus (C\cap B(\tilde x,r))<0\ .$$
This shows that $(b)$ is satisfied. Finally, fix $\alpha, \beta\in L^{\infty}(\Omega)$. Clearly, the function
$$(x,s,t)\to \alpha(x)F(x,s,t)+\beta(x)F(x,s,t)+K(x,s,t)$$ belongs to ${\cal A}$, and so the functional
$$(u,v)\to I(u,v)+\langle \varphi(u,v),(\alpha,\beta)\rangle_Y$$
is sequentially weakly lower semicontinuous in $X$. Let us show that it is coercive. Set
$$\theta=\max\left\{\|\alpha\|_{L^{\infty}(\Omega)},\|\beta\|_{L^{\infty}(\Omega)}\right\}$$
and fix $\epsilon>0$ so that
$$\epsilon<{{1}\over {\theta}}\left ({{\lambda_1}\over {2}}-\eta\right )\ .\eqno{(11)}$$
 By $(a_1)$, there is
$c_{\epsilon}>0$ such that
$$|F(x,s,t)|+|G(x,s,t)|\leq \epsilon(|s|^2+|t|^2)+c_{\epsilon}$$
for all $(x,s,t)\in \Omega\times {\bf R}^2$. Then, for each $u, v\in H^1_0(\Omega)$, recalling $(9)$, we have
$$I(u,v)+\langle \varphi(u,v),(\alpha,\beta)\rangle_Y\geq 
\left ( {{1}\over {2}}-{{\eta}\over {\lambda_1}}\right)\int_{\Omega}(|\nabla u(x)|^2+|\nabla v(x)|^2)dx-
\int_{\Omega}|\alpha(x)F(x,u(x),v(x))+\beta(x)G(x,u(x),v(x))|dx$$
$$\geq \left ( {{1}\over {2}}-{{\eta}\over {\lambda_1}}\right)\int_{\Omega}(|\nabla u(x)|^2+|\nabla v(x)|^2)dx-
\theta\epsilon\int_{\Omega}(|u(x)|^2+|v(x)|^2)dx-\theta c_{\epsilon}\hbox {\rm meas}(\Omega)$$
$$\geq \left ( {{1}\over {2}}-{{\eta}\over {\lambda_1}}-{{\theta\epsilon}\over {\lambda_1}}\right )\int_{\Omega}(|\nabla u(x)|^2+|\nabla v(x)|^2)dx-\theta c_{\epsilon}\hbox {\rm meas}(\Omega)\ .$$
Notice that, in view of $(11)$, we have ${{1}\over {2}}-{{\eta}\over {\lambda_1}}-{{\theta\epsilon}\over {\lambda_1}}>0$, and so
$$\lim_{\|(u,v)\|_X\to +\infty}(I(u,v)+\langle \varphi(u,v),(\alpha,\beta)\rangle_Y)=+\infty\ ,$$
as claimed. In particular, this also implies that the functional $(u,v)\to I(u,v)+\langle \varphi(u,v),(\alpha,\beta)\rangle_Y$ is
weakly lower semicontinuous, by the Eberlein-Smulyan theorem. Thus, the assumptions of Theorem 2.1 are satisfied. Therefore, for each
convex set $S\subseteq L^{\infty}(\Omega)\times L^{\infty}(\Omega)$ dense in $H^1_0(\Omega)\times H^1_0(\Omega)$, there exists
$(\alpha,\beta)\in S$, such that the functional
$$(u,v)\to {{1}\over {2}}\left ( \int_{\Omega}|\nabla u(x)|^2dx+\int_{\Omega}|\nabla v(x)|^2dx\right )-
\int_{\Omega}(\alpha(x)F(x,u(x),v(x))+\beta(x)G(x,u(x),v(x))+K(x,u(x),v(x)))dx$$
has at least two global minima in $H^1_0(\Omega)\times H^1_0(\Omega)$. Finally, by Example 38.25 of [17], the same functional satisfies
the Palais-Smale condition, and so it admits at least three critical points, in view of Corollary 1 of [3]. The proof is complete.\hfill 
$\bigtriangleup$\par
\medskip
REMARK 2.1. - We are not aware of known results close enough to Theorem 2.2 in order to do a proper comparison. This sentence also
applies to the case of single equations, that is to say when $F, G, K$ depend on $x$ and $s$ only. For an account on elliptic systems, we refer
to [2].
\medskip
Among the various corollaries of Theorem 2.2, we wish to stress the following ones:\par
\medskip
COROLLARY 2.1. - {\it Let $K\in {\cal A}$, with $K(x,0,0)=0$ for all $x\in \Omega$, satisfy condition $(a_2)$. Moreover, let
$\Phi:{\bf R}^2\to {\bf R}$ be a non-constant $C^1$ function, with $\Phi(0,0)=0$, belonging to ${\cal A}$, with $p={{2}\over {n-2}}$ when $n>2$.
\par
Then, for every convex set $S\subseteq L^{\infty}(\Omega)\times L^{\infty}(\Omega)$ dense in $L^2(\Omega)\times
L^2(\Omega)$, there exists $(\alpha,\beta)\in S$ such that the problem
 $$\cases {-\Delta u=(\alpha(x)\cos(\Phi(u,v))-\beta(x)\sin(\Phi(u,v)))\Phi_u(u,v)+K_u(x,u,v) & in $\Omega$
\cr & \cr
-\Delta v= (\alpha(x)\cos(\Phi(u,v))-\beta(x)\sin(\Phi(u,v)))\Phi_v(u,v)+K_v(x,u,v) & in $\Omega$
\cr & \cr
u=v=0 & on $\partial\Omega$\cr}$$
has at least three weak solutions, two of which are global minima in $H^1_0(\Omega)\times H^1_0(\Omega)$ of the functional
$$(u,v)\to {{1}\over {2}}\left ( \int_{\Omega}|\nabla u(x)|^2dx+\int_{\Omega}|\nabla v(x)|^2dx\right )$$
$$-\int_{\Omega}(\alpha(x)\sin(\Phi(u(x),v(x)))+\beta(x)\cos(\Phi(u(x),v(x)))+K(x,u(x),v(x)))dx\ .$$}\par
\smallskip
PROOF. It suffices to apply Theorem 2.2 to the functions $F, G:{\bf R}^2\to {\bf R}$ defined by
$$F(s,t)=\sin(\Phi(s,t))\ ,$$
$$G(s,t)=\cos(\Phi(s,t))$$
for all $(s,t)\in {\bf R}^2$.\hfill $\bigtriangleup$
\par
\medskip
COROLLARY 2.2. - {\it Let $F, G:{\bf R}\to {\bf R}$ belong to ${\cal A}$, with $p={{2}\over {n-2}}$ when $n>2$.
Moreover, assume that $F, G$ are twice differentiable at $0$ and that
$$\lim_{|s|\to +\infty}{{|F(s)|+|G(s)|}\over {s^2}}=0 \ ,$$
$$0<|F(0)|^2+|G(0)|^2=\inf_{s\in {\bf R}}(|F(s)|^2+|G(s)|^2)\ ,$$
$$F''(0)F(0)+G''(0)G(0)<0\ .\eqno{(12)}$$
Then, for every convex set $S\subseteq L^{\infty}(\Omega)\times L^{\infty}(\Omega)$ dense in $L^2(\Omega)\times
L^2(\Omega)$, there exists $(\alpha,\beta)\in S$ such that the problem
$$\cases {-\Delta u=\alpha(x)F'(u)+\beta(x)G'(u) & in $\Omega$
\cr & \cr
u=0 & on $\partial\Omega$\cr}$$
has at least three weak solutions, two of which are global minima in $H^1_0(\Omega)$ of the functional
$$u\to {{1}\over {2}}\int_{\Omega}|\nabla u(x)|^2dx-
\int_{\Omega}(\alpha(x)F(u(x))+\beta(x)G(u(x)))dx\ .$$}\par
\smallskip
PROOF. We apply Theorem 2.2 taking $K=0$. 
Since $0$ is a global minimum of the function $|F(\cdot)|^2+|G(\cdot)|^2$, we
have
$$F'(0)F(0)+G'(0)G(0)=0$$
and so, in view of $(12)$, $0$ is a strict local maximum for the function $F(\cdot)F(0)+G(\cdot)G(0)$. Hence, $(a_4)$ is satisfied and Theorem 2.2
gives the conclusion.\hfill $\bigtriangleup$\par
\medskip

\vfill\eject
\centerline {\bf References}\par
\bigskip
\bigskip
\noindent
[1]\hskip 5pt A. CABADA and A. IANNIZZOTTO, {\it A note on a question of Ricceri},
Appl. Math. Lett., {\bf 25} (2012), 215-219.\par
\smallskip
\noindent
[2]\hskip 5pt D. G. de FIGUEIREDO, {\it Semilinear elliptic systems: existence, multiplicity, symmetry of solutions}. Handbook of differential equations: stationary partial differential equations. Vol. V, 1-48, Handb. Differ. Equ., Elsevier/North-Holland, Amsterdam, 2008. \par
\smallskip
\noindent
[3]\hskip 5pt P. PUCCI and J. SERRIN, {\it A mountain pass
theorem}, J. Differential Equations, {\bf 60} (1985), 142-149.\par
\smallskip
\noindent
[4]\hskip 5pt B. RICCERI, {\it On a three critical points theorem}, Arch. Math.,  {\bf 75} (2000),  220-226.\par
\smallskip
\noindent
[5]\hskip 5pt B. RICCERI, {\it Well-posedness of constrained minimization
problems via saddle-points}, J. Global Optim., {\bf 40} (2008),
389-397.\par
\smallskip
\noindent
[6]\hskip 5pt B. RICCERI, {\it Multiplicity of global minima for
parametrized functions}, Rend. Lincei Mat. Appl., {\bf 21} (2010),
47-57.\par
\smallskip
\noindent
[7]\hskip 5pt B. RICCERI, {\it Nonlinear eigenvalue problems},  
in ``Handbook of Nonconvex Analysis and Applications'' 
D. Y. Gao and D. Motreanu eds., 543-595, International Press, 2010.\par
\smallskip
\noindent
[8]\hskip 5pt B. RICCERI, {\it A strict minimax inequality criterion and some of its consequences}, Positivity, {\bf 16} (2012), 455-470.\par
\smallskip
\noindent
[9]\hskip 5pt B. RICCERI, {\it A range property related to non-expansive
operators}, Mathematika, {\bf 60} (2014), 232-236.\par
\smallskip
\noindent
[10]\hskip 5pt B. RICCERI, {\it Singular points of non-monotone potential operators}, J. Nonlinear
Convex Anal., {\bf 16} (2015), 1123-1129.\par
\smallskip
\noindent
[11]\hskip 5pt B. RICCERI, {\it Energy functionals of Kirchhoff-type problems having multiple global minima}, Nonlinear Anal., {\bf 115} (2015),  130-136.\par
\smallskip
\noindent
[12]\hskip 5pt B. RICCERI, {\it On a minimax theorem: an improvement, a new proof and an overview of its applications},
Minimax Theory Appl., {\bf 2} (2017), 99-152.\par
\smallskip
\noindent
[13]\hskip 5pt B. RICCERI, {\it Another multiplicity result for the periodic solutions of certain systems}, Linear Nonlinear Anal., {\bf 5} (2019),
371-378.\par
\smallskip
\noindent
[14]\hskip 5pt B. RICCERI, {\it Miscellaneous applications of certain minimax theorems II}, Acta Math. Vietnam., {\bf 45} (2020), 515-524.\par
\smallskip
\noindent
[15]\hskip 5pt B. RICCERI, {\it An alternative theorem for gradient systems}, Pure Appl. Funct. Anal., to appear.\par
\smallskip
\noindent
[16]\hskip 5pt B. RICCERI, {\it A remark on variational inequalities in small balls},
J. Nonlinear Var. Anal., {\bf 4} (2020), 21-26.\par
\smallskip
\noindent
[17]\hskip 5pt E. ZEIDLER, {\it Nonlinear functional analysis and its
applications}, vol. III, Springer-Verlag, 1985.\par
\smallskip
\noindent

\bigskip
\bigskip
\bigskip
\bigskip
Department of Mathematics and Informatics\par
University of Catania\par
Viale A. Doria 6\par
95125 Catania, Italy\par
{\it e-mail address}: ricceri@dmi.unict.it

\bye